\documentclass[3p,preprint,12pt]{elsarticle}
\makeatletter\if@twocolumn\PassOptionsToPackage{switch}{lineno}\else\fi\makeatother

\usepackage{url}

\usepackage{tabulary,xcolor}
\usepackage{amsfonts,amsmath,amssymb}
\usepackage[T1]{fontenc}
\makeatletter
\let\save@ps@pprintTitle\ps@pprintTitle
\def\ps@pprintTitle{\save@ps@pprintTitle\gdef\@oddfoot{\footnotesize\itshape \null\hfill\today}}
\def\hlinewd#1{%
  \noalign{\ifnum0=`}\fi\hrule \@height #1%
  \futurelet\reserved@a\@xhline}

\AtBeginDocument{\ifNAT@numbers \biboptions{sort&compress}\fi}
\makeatother

  \usepackage{subfig} 
  \usepackage{float}
  
  \usepackage{todonotes}
  \usepackage{bm}
  
\usepackage{amssymb}
\usepackage{amsmath}
\usepackage{amsthm}
\usepackage{mathrsfs}

\usepackage{xcolor,pict2e}

\theoremstyle{plain}
\newtheorem{thm}{Theorem}[section]

\newtheorem{rem}[thm]{Remark}

\usepackage{ifluatex}
\ifluatex
\usepackage{fontspec}
\defaultfontfeatures{Ligatures=TeX}
\usepackage[]{unicode-math}
\unimathsetup{math-style=TeX}
\else 
\usepackage[utf8]{inputenc}
\fi 
\ifluatex\else\usepackage{stmaryrd}\fi

\usepackage{url,multirow,morefloats,floatflt,cancel,tfrupee}
\makeatletter

\AtBeginDocument{\@ifpackageloaded{textcomp}{}{\usepackage{textcomp}}}
\makeatother
\usepackage{colortbl}
\usepackage{xcolor}
\usepackage{pifont}
\usepackage[nointegrals]{wasysym}
\makeatletter

\def\mcWidth#1{\csname TY@F#1\endcsname+\tabcolsep}

\def\cAlignHack{\rightskip\@flushglue\leftskip\@flushglue\parindent\z@\parfillskip\z@skip}
\def\rAlignHack{\rightskip\z@skip\leftskip\@flushglue \parindent\z@\parfillskip\z@skip}

\@ifundefined{etal}{}{}

\usepackage{ifxetex}
\ifxetex\else\if@twocolumn\@ifpackageloaded{stfloats}{}{\usepackage{dblfloatfix}}\fi\fi

\AtBeginDocument{
\expandafter\ifx\csname eqalign\endcsname\relax
\def\eqalign#1{\null\vcenter{\def\\{\cr}\openup\jot\m@th
  \ialign{\strut$\displaystyle{##}$\hfil&$\displaystyle{{}##}$\hfil
      \crcr#1\crcr}}\,}
\fi
}

\AtBeginDocument{%
  \@ifpackageloaded{endfloat}%
   {\renewcommand\efloat@iwrite[1]{\immediate\expandafter\protected@write\csname efloat@post#1\endcsname{}}}{\newif\ifefloat@tables}%
}%

\def\BreakURLText#1{\@tfor\brk@tempa:=#1\do{\brk@tempa\hskip0pt}}
\let\lt=<
\let\gt=>
\def\processVert{\ifmmode|\else\textbar\fi}

\@ifundefined{subparagraph}{
\def\subparagraph{\@startsection{paragraph}{5}{2\parindent}{0ex plus 0.1ex minus 0.1ex}%
{0ex}{\normalfont\small\itshape}}%
}{}

\newcommand\role[1]{\unskip}
\newcommand\aucollab[1]{\unskip}
  
\@ifundefined{tsGraphicsScaleX}{\gdef\tsGraphicsScaleX{1}}{}
\@ifundefined{tsGraphicsScaleY}{\gdef\tsGraphicsScaleY{.9}}{}
\def\checkGraphicsWidth{\ifdim\Gin@nat@width>\linewidth
	\tsGraphicsScaleX\linewidth\else\Gin@nat@width\fi}

\def\checkGraphicsHeight{\ifdim\Gin@nat@height>.9\textheight
	\tsGraphicsScaleY\textheight\else\Gin@nat@height\fi}

\def\fixFloatSize#1{}
\let\ts@includegraphics\includegraphics

\def\inlinegraphic[#1]#2{{\edef\@tempa{#1}\edef\baseline@shift{\ifx\@tempa\@empty0\else#1\fi}\edef\tempZ{\the\numexpr(\numexpr(\baseline@shift*\f@size/100))}\protect\raisebox{\tempZ pt}{\ts@includegraphics{#2}}}}

\AtBeginDocument{\def\includegraphics{\@ifnextchar[{\ts@includegraphics}{\ts@includegraphics[width=\checkGraphicsWidth,height=\checkGraphicsHeight,keepaspectratio]}}}

\DeclareMathAlphabet{\mathpzc}{OT1}{pzc}{m}{it}

\def\URL#1#2{\@ifundefined{href}{#2}{\href{#1}{#2}}}

\def\UrlOrds{\do\*\do\-\do\~\do\'\do\"\do\-}%
\g@addto@macro{\UrlBreaks}{\UrlOrds}

\edef\fntEncoding{\f@encoding}

\makeatother

\newif\ifmultipleabstract\multipleabstractfalse%
\begin{document}

\begin{frontmatter}

\title{A stochastic Hamiltonian formulation applied to dissipative particle dynamics}
    
\author[]{Linyu Peng\corref{c-3c9be7c61ad9}}
\ead{l.peng@mech.keio.ac.jp}\cortext[c-3c9be7c61ad9]{Corresponding author.}
\author[]{Noriyoshi Arai}
\ead{arai@mech.keio.ac.jp}
\author[]{Kenji Yasuoka}
\ead{yasuoka@mech.keio.ac.jp}
    
\address{Department of Mechanical Engineering\unskip, 
    Keio University\unskip, Yokohama 223-8522\unskip, Japan}

\begin{abstract}
In this paper, a stochastic Hamiltonian formulation (SHF) is proposed and applied to dissipative particle dynamics  (DPD) simulations. As an extension of Hamiltonian dynamics to stochastic dissipative systems, the SHF provides necessary foundations and great convenience  for constructing efficient numerical integrators. As a first attempt,  we develop the St\"ormer--Verlet type of schemes based on the SHF, which are structure-preserving for deterministic Hamiltonian systems without external forces, the dissipative forces in DPD. Long-time behaviour of the schemes is shown numerically by studying the damped Kubo oscillator. In particular, the proposed schemes include the conventional  Groot--Warren's modified velocity-Verlet method and a modified version of Gibson--Chen--Chynoweth as special cases. The schemes are applied to DPD simulations and analysed numerically.
\vspace{0.2cm}

 {\bf Keywords:} Dissipative particle dynamics; Hamiltonian mechanics; Stochastic differential equations; St\"ormer--Verlet methods
\end{abstract}

\end{frontmatter}
   
\section{Introduction}

A dissipative particle dynamics (DPD) simulation \cite{Hoogerbrugge1992,Groot1997} method is a type of coarse-grained molecular simulation method, which has proven to be a powerful tool for investigating fluid events occurring on a wide range of spatio-temporal scales compared to all-atom simulations. 
Using DPD method, many studies have been conducted for both the statics and dynamics of complex system at the mesoscopic level, such as unique self-assembled structures formed by nanoparticles or polymers \cite{Chen2020,Huang2020,Arai2020,Cheng2021}, mechanical or rheological properties of soft materials \cite{Nikolov2020,Pan2020,Kobayashi2021}, medical materials and biological functions \cite{Papageorgiou2018,Chen2019,Arai2021,Sicard2021}, and so forth.
Huang $et~al.$ \cite{Huang2020} proposed a method to fabricate various two-dimensional nanostructures using self-assembly of block copolymers and demonstrated it in DPD simulations. 
The simulations showed that surface patterns of three-dimensional nanostructures could be evolved to solve problems in lithography and transistors.
In order to overcome the problem of low toughness in the use of humanoid robotic hands, Pan $et~al.$ \cite{Pan2020}  developed an ultra-tough electric tendon based on spider silk toughened with single-wall carbon nanotubes (SWCNTs). 
In that study, DPD simulations were performed to understand how SWCNTs improve the mechanical properties of the fibers at a molecular level. 
Sicard and Toro-Mendoza \cite{Sicard2021} reported on the computational design of soft nanocarriers using pickering emulsions (nanoparticle armored droplet), able to selectively encapsulate or release a probe load under specific flow conditions. 
They described in detail the mechanisms at play in the formation of pocket-like structures and their stability under external flow. Moreover, the rheological properties of the designed nanocarriers were compared with those of delivery systems used in pharmaceutical and cosmetic technologies.

On the other hand, during the last decades, a lot of efforts have been made for proposing  efficient simulation methods for DPD to achieve simultaneous temperature control and  momentum preservation. Examples include Groot--Warren's modified velocity-Verlet (GW) method \cite{Groot1997}, the method of Gibson--Chen--Chynoweth (GCC) \cite{GCC1999},  and splitting methods \cite{Sha2003,Sha2021}; a review and comparison of commonly used methods for DPD are available in \cite{LS2015}. In the current study, we will show that various velocity-Verlet methods for DPD, including GW and GCC methods, are actually special cases of the St\"ormer--Verlet (SV) schemes for a novel stochastic Hamiltonian formulation (SHF) with dissipative forces which are often called external forces in classical Hamiltonian mechanics; in DPD, these dissipative forces are in fact internal forces (see Section \ref{sec:SHF}). To be consistent, they will be called external forces in the general setting but dissipative forces in DPD.
SV schemes are well-known symplectic-preserving numerical methods for deterministic Hamiltonian systems without external forces. 

Symplecticity is a crucial  feature  of Hamiltonian systems. Geometrically, it implies area or volume preservation of the corresponding phase flows due to Liouville's Theorem. Symplectic integrators are among the most important types of geometric numerical integrators for Hamiltonian systems \cite{LR2004,HLW2006}.  Symplectic integrators for stochastic Hamiltonian systems with or without external forces have received great attention as well, e.g., \cite{MRT2002,Wan2009,HT2018,KT2021,MW2001}.  The SHF we propose in the current study can be viewed as a matrix generalisation of stochastic forced Hamiltonian systems studied in \cite{KT2021}; see also \cite{HT2018,LO2008}.    The Hamiltonian structure brings us a convenient setting  for  analysis of the underlying dynamical system; moreover, it allows the systematic construction of structure-preserving  integrators possible. In this paper, we will mainly be  focused on the extension of SV type of symplectic schemes to systems of SHF and to DPD. 

The paper is organised as follows. In Section \ref{sec:SHF}, we propose the SHF  and derive the DPD by specifying the Hamiltonian functions and external/dissipative forces properly. SV type of schemes for the SHF and the DPD are constructed in Section \ref{sec:SV} and in particular, we will be focused on several explicit schemes that are applied to  DPD simulations in Section \ref{sec:app}. Finally, we conclude and point out some future researches  in Section \ref{sec:con}.

\section{The stochastic Hamiltonian formulation with external forces}
\label{sec:SHF}

Let $Q$ be an $n$-dimensional  configuration space of a mechanical system with $\bm{q}$ the generalised coordinates. Let  $(\bm{q},\dot{\bm{q}})\in TQ$ and $(\bm{q},\bm{p})\in T^*Q$ be coordinates of the tangent bundle and the cotangent bundle, respectively. We propose a stochastic Hamiltonian formulation (SHF) with external forces  as a dynamical system in  $T^*Q$ as follows:
\begin{equation}\label{eq:SHF}
\begin{aligned}
\left(\begin{array}{c}
\operatorname{d}\!{\bm{q}}\\ \operatorname{d}\!{\bm{p}}
\end{array}\right)  =J\nabla H(\bm{q},\bm{p}) \operatorname{d}\!t&+\left(\begin{array}{c}
0\\
\bm{F}^{\operatorname{D}}(\bm{q},\bm{p})
\end{array}\right)\operatorname{d}\!t\\
& +\sum_{i=1}^K\sum_{j=1}^K\left(J\nabla h_{ij}(\bm{q},\bm{p})+\left(\begin{array}{c}
0\\
\bm{F}^{\operatorname{SD}}_{ij}(\bm{q},\bm{p})
\end{array}\right)\right)\circ {\operatorname{d}\! W_{ij}(t)},
\end{aligned}
\end{equation}
where $\circ$ denotes the Stratonovich integration, $J$ is the canonical symplectic matrix
\begin{equation}
J=\left(\begin{array}{cc}
0 & I_n\\
-I_n & 0
\end{array}\right),
\end{equation}
$\bm{F}^{\operatorname{D}}:T^*Q\rightarrow T^*Q$ and $\bm{F}_{ij}^{\operatorname{SD}}:T^*Q\rightarrow T^*Q$ are fibre-preserving maps of the external forces leading to dissipation,
the functions $H:T^*Q\rightarrow \mathbb{R}$ and $h_{ij}:T^*Q\rightarrow \mathbb{R}$ are the Hamiltonian functions, and  components of the symmetric $K\times K$ random matrix $W(t)$ are independent Wiener processes. Note that the indices $i,j$ are not necessary of the same dimension.The superindices $\operatorname{D}$ and $\operatorname{SD}$ are shorthand for `Dissipation' and `Stochastic Dissipation', respectively.  For more details on stochastic differential equations, the reader may refer to \cite{Arnold1994,Evans2013,KP1999}.

The SHF \eqref{eq:SHF} can be written in the It\^{o} form as 
\begin{equation}\label{eq:ito}
\operatorname{d}\!\bm{z}=A(\bm{z})\operatorname{d}\!t+\sum_{i=1}^K\sum_{j=1}^KB_{ij}(\bm{z})\operatorname{d}\!W_{ij}(t),
\end{equation}
where $\bm{z}=(\bm{q},\bm{p})^{\operatorname{T}}$, 
\begin{equation}\label{eq:AA}
A(\bm{z})=\left(\begin{array}{c}
\nabla_{\bm{p}}H+\frac{1}{2}\sum\limits_{i=1}^K\sum\limits_{j=1}^K\left(\frac{\partial^2 h_{ij}}{\partial \bm{p}\partial\bm{q}}\left(\nabla_{\bm{p}}h_{ij}\right)+\frac{\partial^2 h_{ij}}{\partial\bm{p}^2}\left(\bm{F}_{ij}^{\operatorname{SD}}-\nabla_{\bm{q}}h_{ij}\right)\right)\\
\hspace{-1cm} -\nabla_{\bm{q}}H+\bm{F}^{\operatorname{D}}+\frac{1}{2}\sum\limits_{i=1}^K\sum\limits_{j=1}^K\left(\left(\frac{\partial^2 h_{ij}}{\partial \bm{q}\partial\bm{p}}-\nabla_{\bm{\bm{p}}}\bm{F}_{ij}^{\operatorname{SD}}\right)\left(\nabla_{\bm{p}}h_{ij}-\bm{F}_{ij}^{\operatorname{SD}}\right)\right.\\
 \quad \quad\quad \quad\quad \quad\quad \quad\quad \quad\quad \quad\quad \quad\quad \quad \left.-\left(\frac{\partial^2 h_{ij}}{\partial\bm{q}^2}-\nabla_{\bm{q}}\bm{F}_{ij}^{\operatorname{SD}}\right)\left(\nabla_{\bm{p}}h_{ij}\right)\right)
\end{array}
\right)
\end{equation}
and 
\begin{equation}
B_{ij}(\bm{z})=\left(\begin{array}{c}
\nabla_{\bm{p}}h_{ij}\\
-\nabla_{\bm{q}}h_{ij}+\bm{F}_{ij}^{\operatorname{SD}}
\end{array}
\right).
\end{equation}
Here, ${\partial^2 h_{ij}}/{\partial \bm{p}\partial\bm{q}}$,   ${\partial^2 h_{ij}}/{\partial\bm{q}}^2$ and  ${\partial^2 h_{ij}}/{\partial \bm{p}}^2$ denote the Hessian matrices of $h_{ij}$, and $\nabla$ denotes the gradient of functions. Throughout the paper, we will employ the conventional assumptions that the Hamiltonians $H$
 and $h_{ij}$ are all $C^2$ functions and $A$ and $B_{ij}$ are globally Lipschitz \cite{Arnold1994,KP1999}.

\begin{rem}
The SHF can be derived through variational calculus. It will be called a stochastic Lagrange--d'Alembert principle in the phase space $T^*Q$, reading
\begin{equation}\label{eq:vp}
\begin{aligned}
\delta & \int_{t_a}^{t_b} \left(\bm{p}\circ \operatorname{d}\!{\bm{q}} -H(\bm{q},\bm{p}) \operatorname{d}\!t\right)+\int_{t_a}^{t_b} \bm{F}^{\operatorname{D}}(\bm{q},\bm{p})\cdot \delta \bm{q}\operatorname{d}\!t\\
&~~+ \sum_{i=1}^K\sum_{j=1}^K\left(\delta\int_{t_a}^{t_b}-h_{ij}(\bm{q},\bm{p})\circ \operatorname{d}\!W_{ij}(t)+\int_{t_a}^{t_b}\left(\bm{F}_{ij}^{\operatorname{SD}}(\bm{q},\bm{p})\cdot \delta \bm{q}\right)\circ \operatorname{d}\! W_{ij}(t)\right)=0.
\end{aligned}
\end{equation}
The time interval is $[t_a,t_b]$ ($t_a<t_b$).
The first row denotes all deterministic terms, while the second row includes all stochastic terms. 

Solutions of the SHF \eqref{eq:SHF} satisfies the stochastic Lagrange--d'Alembert principle \eqref{eq:vp}; see, e.g., \cite{KT2021}. The converse is also true providing the regularity of $\bm{q}$ and $\bm{p}$  \cite{SC2021}. In particular if $h_{ij}=h_{ij}(\bm{q})$ are all independent of $\bm{p}$, which is exactly the case for DPD, $(\bm{q},\bm{p})$ is a solution of the SHF \eqref{eq:SHF} if and only if it satisfies the stochastic Lagrange--d'Alembert principle \eqref{eq:vp} \cite{Tyr2021}.
 
\end{rem}

{\bf DPD derived from the SHF.}
To derive the DPD system of $N$ particles, we assume that there exist no stochastic dissipative forces, meaning that 
\begin{equation}
\bm{F}^{\operatorname{SD}}_{ij}(\bm{q},\bm{p})\equiv 0, \quad \forall i,j=1,2,\ldots,N.
\end{equation}
 In the general SHF formulation  \eqref{eq:SHF}, introduce the local coordinates for the cotangent bundle of $N$ copies of $Q$ as
\begin{equation}
\bm{q}=(\bm{q}_1,\bm{q}_2,\ldots,\bm{q}_N),\quad \bm{p}=(\bm{p}_1,\bm{p}_2,\ldots,\bm{p}_N),
\end{equation}
where $(\bm{q}_i,\bm{p}_i)$ are
the coordinates of the phase space $T^*Q$ for the $i$-th particle. As commonly considered in DPD, we will  be focused on the three-dimensional Euclidean space,  i.e., $Q=\mathbb{R}^3$, in the current study. Define the Hamiltonian $H(\bm{q},\bm{p})$ as the total energy:
\begin{equation}\label{eq:H}
H(\bm{q},\bm{p})=\sum_{i=1}^N\frac{1}{2m_i}|\bm{p}_i|^2+V(\bm{q}),
\end{equation}
where the potential energy $V(\bm{q})$ is given by
\begin{equation}
V(\bm{q})=\sum_{i=1}^N\sum_{j=1}^N\frac{a_{ij}}{4}q_{\mathrm{c}}\left(1-\frac{q_{ij}}{q_{\mathrm{c}}}\right)^2\delta_{ij}.
\end{equation}
Here  $m_i$ is mass of the $i$-th particle, $q_{\mathrm{c}}$ is a constant, $a_{N\times N}$ is a constant symmetric matrix, $q_{ij}=|\bm{q}_i-\bm{q}_j|$ is the distance of the $i$-th and the $j$-th particles, and $\delta_{ij}$ is given by
\begin{equation}
\delta_{ij}=\left\{\begin{array}{cl}
1, & q_{ij}<q_{\mathrm{c}}, \vspace{0.2cm}\\
 0, & q_{ij}\geq q_{\mathrm{c}}.
\end{array}\right.
\end{equation}

\begin{rem} Direct computation gives  gradient of the Hamiltonian $H$ as follows
\begin{equation}
\nabla H(\bm{q},\bm{p})=\left(-\sum_{j\neq i} \bm{F}_{ij}^{\operatorname{C}}(\bm{q}),\frac{\bm{p}_i}{m_i}\right)^{\operatorname{T}},
\end{equation}
where the conservative force reads
\begin{equation}\label{eq:FC}
\bm{F}^{\operatorname{C}}_{ij}(\bm{q})=a_{ij}\left(1-\frac{q_{ij}}{q_{\mathrm{c}}}\right)\delta_{ij}\bm{\widehat{q}}_{ij},\quad i,j=1,2,\ldots, N,\quad i\neq j,
\end{equation}
with 
\begin{equation}
 \bm{\widehat{q}}_{ij}=\frac{\bm{q}_i-\bm{q}_j}{q_{ij}}=\frac{\bm{q}_i-\bm{q}_j}{|\bm{q}_i-\bm{q}_j|}
 \end{equation}
 and the superindex $\operatorname{C}$ meaning `Conservation'. Obviously,  the conservative force $\bm{F}^{\operatorname{C}}_{ij}(\bm{q})$ between the $i$-th and the $j$-th particles only depends on their relative distance $\bm{q}_i-\bm{q}_j$.
\end{rem}

Furthermore, the (deterministic) dissipative force is defined by
\begin{equation}\label{eq:FD}
\bm{F}^{\operatorname{D}}_i(\bm{q},\bm{p})=-\gamma \sum_{j\neq i}\omega^{\operatorname{D}}(q_{ij})\left(\bm{\widehat{q}}_{ij}\cdot \bm{v}_{ij}\right)\bm{\widehat{q}}_{ij},\quad i=1,2,\ldots,N,
\end{equation}
where $\gamma$ is a constant friction parameter,
\begin{equation}
 \bm{v}_{ij}=\frac{\bm{p}_i}{m_i}-\frac{\bm{p}_j}{m_j}
\end{equation}
and  $\omega^{\operatorname{D}}(q_{ij})=\left(\omega^{\operatorname{R}}(q_{ij})\right)^2$ with
\begin{equation}
\omega^{\operatorname{R}}(q_{ij})=\left(1-\frac{q_{ij}}{q_{\mathrm{c}}}\right)\delta_{ij}.
\end{equation}
Here, the superindex $\operatorname{R}$ means `Randomness'.

Let $k=N$ and define the Hamiltonian functions $h_{ij}(\bm{q},\bm{p})$ ($i,j,=1,2,\ldots, N$)  by
\begin{equation}\label{eq:hij}
h_{ij}(\bm{q})=\frac{\sigma}{4}q_{\mathrm{c}}\left(1-\frac{q_{ij}}{q_{\mathrm{c}}}\right)^2\delta_{ij}, 
\end{equation}
where $\sigma$ is a constant noise parameter. Obviously, $h_{ii}(\bm{q})\equiv \text{const}$ for all $i=1,2,\ldots,N$ and hence $\nabla h_{ii}(\bm{q})\equiv 0$.

\begin{rem}
When $i\neq j$, since the Hamiltonian function $h_{ij}(\bm{q})$ only depends on $\bm{q}_i$ and $\bm{q}_j$, nonzero components of  its gradient are given by
\begin{equation}
\begin{aligned}
\nabla_{\bm{q}_i} h_{ij}(\bm{q})&=-\frac{\sigma}{2}\left(1-\frac{q_{ij}}{q_{\mathrm{c}}}\right)\delta_{ij}\bm{\widehat{q}}_{ij}=-\frac{\sigma}{2}\omega^{\operatorname{R}}(q_{ij})\bm{\widehat{q}}_{ij},\\
\nabla_{\bm{q}_j} h_{ij}(\bm{q})&=-\frac{\sigma}{2}\left(1-\frac{q_{ij}}{q_{\mathrm{c}}}\right)\delta_{ij}\bm{\widehat{q}}_{ji}=-\frac{\sigma}{2}\omega^{\operatorname{R}}(q_{ij})\bm{\widehat{q}}_{ji}.
\end{aligned}
\end{equation}
\end{rem}


Substituting the functions specified above to the SHF \eqref{eq:SHF}, we obtain the system of DPD as follows
\begin{equation}\label{eq:DPD}
\left\{ \begin{aligned}
\dot{\bm{q}}_i&=\frac{\bm{p}_i}{m_i},\\
\dot{\bm{p}}_i&=\sum_{j\neq i}\bm{F}_{ij}^{\operatorname{C}}(\bm{q})+\bm{F}^{\operatorname{D}}_i(\bm{q},\bm{p})+\sigma\sum_{j\neq i}\omega^{\operatorname{R}}(q_{ij})\bm{\widehat{q}}_{ij}\circ \frac{\operatorname{d}\!W_{ij}(t)}{\operatorname{d}\!t},
\end{aligned}
 \right.
\end{equation}
for $i=1,2,\ldots,N$, in which the dissipative   force $\bm{F}^{\operatorname{D}}_i(\bm{q},\bm{p})$ is given by \eqref{eq:FD} while the  conservative force $\bm{F}^{\operatorname{C}}_i(\bm{q},\bm{p})$ and randomness contribution are respectively derived from the Hamiltonians $H(\bm{q},\bm{p})$, i.e., the total energy \eqref{eq:H}, and $h_{ij}(\bm{q},\bm{p})$ defined in \eqref{eq:hij}.   It is obvious that the DPD system \eqref{eq:DPD} can also be obtained via the stochastic Lagrange--d'Alembert principle \eqref{eq:vp}.  Note that the SHF \eqref{eq:SHF} is {\it formally} divided by $\operatorname{d}\!t$ on both sides to obtain the system \eqref{eq:DPD}, which has been the conventional form of DPD. 

\begin{rem}
Since no stochastic dissipative forces exist and $h_{ij}=h_{ij}(\bm{q})$ are independent from $\bm{p}$, SHF's  It\^o form \eqref{eq:ito}, in particular the coefficient matrix $A(\bm{z})$ given by \eqref{eq:AA},  yields that the DPD \eqref{eq:DPD} takes the same form in both the It\^o framework and the Stratonovich framework.
\end{rem}



\section{St\"ormer--Verlet schemes for the SHF and the DPD}
\label{sec:SV}

In this section, we propose the St\"ormer--Verlet (SV) type of symplectic schemes for the DPD based on the SHF \eqref{eq:SHF}. That is, when no external forces and randomness are imposed, the corresponding discrete `flow' shall be symplectic as well. In other words, dissipation in the numerical schemes is only contributed by the external forces,  same as what occurs in the continuous counterpart.

\subsection{SV type of schemes for the SHF}

Discretize the time interval $[t_a,t_b]$ as a series $t_a=t_0,t_1,t_2,\ldots,t_K=t_b$ and denote $$\Delta t=t_{k+1}-t_k=\frac{t_b-t_a}{K}$$ as the time step. The space $TT^*Q$ where SHF systems (and the corresponding variational structure) are defined is discretized into two copies of the cotangent bundle, i.e., $T^*Q\times T^*Q$, with local coordinates $(\bm{q}^k,\bm{p}^k,\bm{q}^{k+1},\bm{p}^{k+1})$ where $\bm{q}^k=\bm{q}(t_k)$, $\bm{p}^k=\bm{p}(t_k)$ and so forth. 
In the current paper, we will mainly be focused on extensions of  the SV schemes for systems of SHF \eqref{eq:SHF}, which are symplectic schemes of second order accuracy for conservative Hamiltonian systems. The SV schemes arise as the composite of Euler-A and Euler-B methods which are both symplectic, implicit and of first order accuracy for conservative Hamiltonian systems. We will follow a similar approach to introduce SV schemes for the SHF.

For SHF \eqref{eq:SHF}, we propose a family
of  Euler-A methods:
\begin{equation}\label{alg:EA}
\begin{aligned}
{\bm{q}^{k+1}-\bm{q}^k}&=\Delta t * \nabla_{\bm{p}} H(\bm{q}^{k+1},\bm{p}^k)+\sum_{i,j}\nabla_{\bm{p}}h_{ij}(\bm{q}^{k+1},\bm{p}^k)\circ \Delta W_{ij}(t_k),\\
{\bm{p}^{k+1}-\bm{p}^k}&=- \Delta t \left( \nabla_{\bm{q}}H(\bm{q}^{k+1},\bm{p}^k)+(\bm{F}^{\operatorname{D}}(\bm{q},\bm{p}))^k\right)
\\&~~~~+\sum_{i,j}\left(-\nabla_{\bm{q}}h_{ij}(\bm{q}^{k+1},\bm{p}^k)+(\bm{F}_{ij}^{\operatorname{SD}}(\bm{q},\bm{p}))^k\right)\circ \Delta W_{ij}(t_k),
\end{aligned}
\end{equation}
%
%
%
and a family of Euler-B methods:
\begin{equation}\label{alg:EB}
\begin{aligned}
{\bm{q}^{k+1}-\bm{q}^k}&=\Delta t* \nabla_{\bm{p}} H(\bm{q}^{k},\bm{p}^{k+1})+\sum_{i,j}\nabla_{\bm{p}}h_{ij}(\bm{q}^{k},\bm{p}^{k+1})\circ \Delta W_{ij}(t_k),\\
{\bm{p}^{k+1}-\bm{p}^k}&=- \Delta t \left(\nabla_{\bm{q}}H(\bm{q}^{k},\bm{p}^{k+1})+(\bm{F}^{\operatorname{D}}(\bm{q},\bm{p}))^k\right)
\\&~~~~+\sum_{i,j}\left(-\nabla_{\bm{q}}h_{ij}(\bm{q}^{k},\bm{p}^{k+1})+(\bm{F}_{ij}^{\operatorname{SD}}(\bm{q},\bm{p}))^k\right)\circ \Delta W_{ij}(t_k),
\end{aligned}
\end{equation}
where $(\bm{F}^{\operatorname{D}}(\bm{q},\bm{p}))^k$ and $(\bm{F}_{ij}^{\operatorname{SD}}(\bm{q},\bm{p}))^k$ denote discretisations of the external forces, and 
\begin{equation}
\Delta W_{ij}(t_k)=W_{ij}(t_{k+1})-W_{ij}(t_k)\sim \mathcal{N}(0,\Delta t).
\end{equation}
Here $\mathcal{N}(0,\Delta t)$ denotes the normal distribution with mean $0$ and standard deviation $\sqrt{\Delta t}$. 

Two types of  SV schemes can be defined as  composites of the Euler methods with time step $\Delta t/2$, namely  $\text{(Euler-A)}\circ \text{(Euler-B)}$ and $\text{(Euler-B)}\circ \text{(Euler-A)}$, which will be called SV-AB schemes and SV-BA schemes, respectively. 

The family of SV-AB schemes, namely $\text{(Euler-A)}\circ \text{(Euler-B)}$, reads
\begin{equation}\label{alg:SVAB}
\begin{aligned}
\bm{p}^{k+1/2}&-\bm{p}^k=\frac{\Delta t}{2}\left[- \nabla_{\bm{q}}H(\bm{q}^{k},\bm{p}^{k+1/2})+(\bm{F}^{\operatorname{D}}(\bm{q},\bm{p}))^{k_1}\right]
\\&\quad \quad \quad +\sum_{i,j}\left(-\nabla_{\bm{q}}h_{ij}(\bm{q}^{k},\bm{p}^{k+1/2})+(\bm{F}_{ij}^{\operatorname{SD}}(\bm{q},\bm{p}))^{k_1}\right)\circ \overline{\Delta} W_{ij}(t_k),\\
\bm{q}^{k+1}&-\bm{q}^k=\frac{\Delta t}{2} \left[ \nabla_{\bm{p}}H(\bm{q}^k,\bm{p}^{k+1/2}) +\nabla_{\bm{p}} H(\bm{q}^{k+1},\bm{p}^{k+1/2}) \right] \\&+\sum_{i,j}\nabla_{\bm{p}}h_{ij}(\bm{q}^k,\bm{p}^{k+1/2})\circ \overline{\Delta} W_{ij}(t_k) +\sum_{i,j}\nabla_{\bm{p}}h_{ij}(\bm{q}^{k+1},\bm{p}^{k+1/2})\circ \overline{\Delta}W_{ij}(t_{k+1/2}),\\
\bm{p}^{k+1}&-\bm{p}^{k+1/2}=\frac{\Delta t}{2}\left[-\nabla_{\bm{q}} H(\bm{q}^{k+1},\bm{p}^{k+1/2})+(\bm{F}^{\operatorname{D}}(\bm{q},\bm{p}))^{k_2}\right]\\
&\quad \quad \quad +\sum_{i,j}\left(-\nabla_{\bm{q}}h_{ij}(\bm{q}^{k+1},\bm{p}^{k+1/2})+(\bm{F}_{ij}^{\operatorname{SD}}(\bm{q},\bm{p}))^{k_2}\right)\circ \overline{\Delta} W_{ij}(t_{k+1/2}),
\end{aligned}
\end{equation}
while the family of SV-BA schemes, namely $\text{(Euler-B)}\circ \text{(Euler-A)}$, reads
\begin{equation}\label{alg:SVBA}
\begin{aligned}
\bm{q}^{k+1/2}-\bm{q}^k&=\frac{\Delta t}{2} * \nabla_{\bm{p}}H(\bm{q}^{k+1/2},\bm{p}^{k})+\sum_{i,j}\nabla_{\bm{p}}h_{ij}(\bm{q}^{k+1/2},\bm{p}^{k})\circ \overline{ \Delta} W_{ij}(t_k),\\
\bm{p}^{k+1}-\bm{p}^k&= \frac{\Delta t}{2} \left[ -\nabla_{\bm{q}}H(\bm{q}^{k+1/2},\bm{p}^{k}) -\nabla_{\bm{q}}H(\bm{q}^{k+1/2},\bm{p}^{k+1}) \right]+ \Delta t * (\bm{F}^{\operatorname{D}}(\bm{q},\bm{p}))^k\\
&~~~~+  \sum_{i,j}\left(-\nabla_{\bm{q}}h_{ij}(\bm{q}^{k+1/2},\bm{p}^{k})+(\bm{F}_{ij}^{\operatorname{SD}}(\bm{q},\bm{p}))^{k_1}\right)\circ \overline{\Delta} W_{ij}(t_k) \\
&~~~~+\sum_{i,j}\left(-\nabla_{\bm{q}}h_{ij}(\bm{q}^{k+1/2},\bm{p}^{k+1})+(\bm{F}_{ij}^{\operatorname{SD}}(\bm{q},\bm{p}))^{k_2}\right)\circ \overline{\Delta}W_{ij}(t_{k+1/2}),\\
\bm{q}^{k+1}-\bm{q}^{k+1/2}&=\frac{\Delta t}{2}*\nabla_{\bm{p}}H(\bm{q}^{k+1/2},\bm{p}^{k+1})+\sum_{i,j}\nabla_{\bm{p}}h_{ij}(\bm{q}^{k+1/2},\bm{p}^{k+1})\circ \overline{ \Delta} W_{ij}(t_{k+1/2}),
\end{aligned}
\end{equation}
where $(\bm{F}^{\operatorname{D}}(\bm{q},\bm{p}))^k$, $(\bm{F}^{\operatorname{D}}(\bm{q},\bm{p}))^{k_1}$ and $(\bm{F}^{\operatorname{D}}(\bm{q},\bm{p}))^{k_2}$ denote three independent discretisations of the force $\bm{F}^{\operatorname{D}}(\bm{q},\bm{p})$,  $(\bm{F}^{\operatorname{SD}}_{ij}(\bm{q},\bm{p}))^{k_1}$ and $(\bm{F}^{\operatorname{SD}}_{ij}(\bm{q},\bm{p}))^{k_2}$ denote two independent discretisations of the force $\bm{F}^{\operatorname{SD}}_{ij}(\bm{q},\bm{p})$, and
\begin{equation}
\overline{\Delta}W_{ij}(t_k) = W_{ij}(t_{k+1/2})-W_{ij}(t_k)\sim \mathcal{N}(0, \Delta t/2).
\end{equation}

\begin{rem}  It is obvious that the SV schemes  \eqref{alg:SVAB} and \eqref{alg:SVBA} reduce to the ordinary SV schemes for conservative Hamiltonian systems, assuming the absence of  external forces and stochastic terms. Consequently, discretisations of the external forces can, in principle, be chosen arbitrarily, providing the resulting schemes are stable and convergent.  Only when discretisations of the external forces are specified properly, they are a 2-stage stochastic partitioned Runge--Kutta method given in \cite{KT2021}; however, in DPD simulations, for instance the GW and GCC methods, these discretisations are often chosen very differently as we will find out below. 

\end{rem}
 
{\bf Separable Hamiltonians.} Assuming the Hamiltonians can be  separated as 
\begin{equation}
H(\bm{q},\bm{p})=T(\bm{p})+V(\bm{q}) \text{ and } h_{ij}(\bm{q},\bm{p})=S_{ij}(\bm{p})+U_{ij}(\bm{q}),
\end{equation}
  the SV-AB schemes \eqref{alg:SVAB} and SV-BA schemes \eqref{alg:SVBA} become
\begin{equation}\label{alg:SVABe}
\begin{aligned}
\bm{p}^{k+1/2}-\bm{p}^k&=\frac{\Delta t}{2}\left[- \nabla_{\bm{q}}V(\bm{q}^{k})+(\bm{F}^{\operatorname{D}}(\bm{q},\bm{p}))^{k_1}\right]
\\&\quad \quad  +\sum_{i,j}\left(-\nabla_{\bm{q}}U_{ij}(\bm{q}^{k})+(\bm{F}_{ij}^{\operatorname{SD}}(\bm{q},\bm{p}))^{k_1}\right)\circ \overline{\Delta} W_{ij}(t_k),\\
\bm{q}^{k+1}-\bm{q}^k&=\Delta t * \nabla_{\bm{p}}T(\bm{p}^{k+1/2}) +\sum_{i,j}\nabla_{\bm{p}}S_{ij}(\bm{p}^{k+1/2})\circ \Delta W_{ij}(t_k),\\
\bm{p}^{k+1}-\bm{p}^{k+1/2}&=\frac{\Delta t}{2}\left[-\nabla_{\bm{q}} V(\bm{q}^{k+1})+(\bm{F}^{\operatorname{D}}(\bm{q},\bm{p}))^{k_2}\right]\\
&\quad \quad  +\sum_{i,j}\left(-\nabla_{\bm{q}}U_{ij}(\bm{q}^{k+1})+(\bm{F}_{ij}^{\operatorname{SD}}(\bm{q},\bm{p}))^{k_2}\right)\circ \overline{\Delta} W_{ij}(t_{k+1/2}),
\end{aligned}
\end{equation}
and
\begin{equation}\label{alg:SVBAe}
\begin{aligned}
\bm{q}^{k+1/2}-\bm{q}^k&=\frac{\Delta t}{2} * \nabla_{\bm{p}}T(\bm{p}^{k})+\sum_{i,j}\nabla_{\bm{p}}S_{ij}(\bm{p}^{k})\circ \overline{ \Delta} W_{ij}(t_k),\\
\bm{p}^{k+1}-\bm{p}^k&=\Delta t * \left( -  \nabla_{\bm{q}}V(\bm{q}^{k+1/2}) + (\bm{F}^{\operatorname{D}}(\bm{q},\bm{p}))^k \right)-\sum_{ij} \nabla_{\bm{q}}U_{ij}(\bm{q}^{k+1/2})\circ \Delta W_{ij}(t_k)\\
&~~~~+  \sum_{i,j}\left((\bm{F}_{ij}^{\operatorname{SD}}(\bm{q},\bm{p}))^{k_1}\circ \overline{\Delta} W_{ij}(t_k)+(\bm{F}_{ij}^{\operatorname{SD}}(\bm{q},\bm{p}))^{k_2}\circ \overline{\Delta}W_{ij}(t_{k+1/2})\right),\\
\bm{q}^{k+1}-\bm{q}^{k+1/2}&=\frac{\Delta t}{2}*\nabla_{\bm{p}}T(\bm{p}^{k+1})+\sum_{i,j}\nabla_{\bm{p}}S_{ij}(\bm{p}^{k+1})\circ \overline{ \Delta} W_{ij}(t_{k+1/2}).
\end{aligned}
\end{equation}




\subsection{SV schemes for the DPD} 
\label{subsec:SVDPD}
The Hamiltonians \eqref{eq:H} and \eqref{eq:hij} of the DPD \eqref{eq:DPD} can obviously be separated with respect to their position and momentum coordinates. Substituting them into  \eqref{alg:SVABe}  and \eqref{alg:SVBAe}, the  SV-AB schemes and SV-BA schemes for DPD turn out to be  
\begin{equation}\label{alg:SV-AB-DPD}
\begin{aligned}
{\bm{p}^{k+1/2}_i-\bm{p}^k_i}&=\frac{\Delta t}{2}\left[ \sum_{j\neq i}\bm{F}_{ij}^{\operatorname{C}}(\bm{q}^{k})+(\bm{F}^{\operatorname{D}}_i(\bm{q},\bm{p}))^{k_1}\right]+\sigma\sum_{j\neq i}\omega^{\operatorname{R}}(q_{ij}^{k})\bm{\widehat{q}}_{ij}^{k}\circ \overline{\Delta} W_{ij}(t_k),\\
{\bm{q}^{k+1}_i-\bm{q}^k_i}&=\Delta t * \frac{\bm{p}_i^{k+1/2}}{m_i},\\
{\bm{p}^{k+1}_i-\bm{p}^{k+1/2}_i}&=\frac{\Delta t}{2}\left[\sum_{j\neq i}\bm{F}_{ij}^{\operatorname{C}}(\bm{q}^{k+1})+(\bm{F}^{\operatorname{D}}_i(\bm{q},\bm{p}))^{k_2}\right]
+\sigma\sum_{j\neq i}\omega^{\operatorname{R}}(q_{ij}^{k+1})\bm{\widehat{q}}_{ij}^{k+1}\circ \overline{\Delta} W_{ij}(t_{k+1/2}),
\end{aligned}
\end{equation}
and
\begin{equation}\label{alg:SV-BA-DPD}
\begin{aligned}
\bm{q}^{k+1/2}_i-\bm{q}^k_i&=\frac{\Delta t}{2} * \frac{\bm{p}_i^k}{m_i},\\
\bm{p}^{k+1}_i-\bm{p}^k_i&=\Delta t \left(  \sum_{j\neq i}\bm{F}_{ij}^{\operatorname{C}}(\bm{q}^{k+1/2})+ (\bm{F}_i^{\operatorname{D}}(\bm{q},\bm{p}))^k \right) +\sigma\sum_{j\neq i}\omega^{\operatorname{R}}(q_{ij}^{k+1/2})\bm{\widehat{q}}_{ij}^{k+1/2}\circ \Delta W_{ij}(t_k),\\
\bm{q}^{k+1}_i-\bm{q}^{k+1/2}_i&=\frac{\Delta t}{2}*\frac{\bm{p}_i^{k+1}}{m_i}.
\end{aligned}
\end{equation}




\begin{rem}\label{rem:MSV}
If we (partially) eliminate the half values in the SV-AB schemes \eqref{alg:SV-AB-DPD} and SV-BA schemes \eqref{alg:SV-BA-DPD}, we can rewrite them in the following equivalent representatives
\begin{equation}\label{MSV-AB}
\begin{aligned}
\bm{q}^{k+1}_i-\bm{q}^k_i&=\frac{\Delta t}{m_i}\left\{\bm{p}_i^k+\frac{\Delta t}{2}\left[ \sum_{j\neq i}\bm{F}_{ij}^{\operatorname{C}}(\bm{q}^{k})+(\bm{F}^{\operatorname{D}}_i(\bm{q},\bm{p}))^{k_1}\right]+\sigma\sum_{j\neq i}\omega^{\operatorname{R}}(q_{ij}^{k})\bm{\widehat{q}}_{ij}^{k}\circ \overline{\Delta} W_{ij}(t_k)\right\},\\
\bm{p}^{k+1}_i-\bm{p}^{k}_i&=\frac{\Delta t}{2}\left[ \sum_{j\neq i}\bm{F}_{ij}^{\operatorname{C}}(\bm{q}^{k})+(\bm{F}^{\operatorname{D}}_i(\bm{q},\bm{p}))^{k_1}\right]+\sigma\sum_{j\neq i}\omega^{\operatorname{R}}(q_{ij}^{k})\bm{\widehat{q}}_{ij}^{k}\circ \overline{\Delta} W_{ij}(t_k)\\
&\quad \quad +\frac{\Delta t}{2}\left[\sum_{j\neq i}\bm{F}_{ij}^{\operatorname{C}}(\bm{q}^{k+1})+(\bm{F}^{\operatorname{D}}_i(\bm{q},\bm{p}))^{k_2}\right]+\sigma\sum_{j\neq i}\omega^{\operatorname{R}}(q_{ij}^{k+1})\bm{\widehat{q}}_{ij}^{k+1}\circ \overline{\Delta} W_{ij}(t_{k+1/2})
\end{aligned}
\end{equation}
and
\begin{equation}\label{MSV-BA}
\begin{aligned}
\bm{q}^{k+1/2}_i-\bm{q}^k_i&=\frac{\Delta t}{2} * \frac{\bm{p}_i^k}{m_i},\\
\bm{p}^{k+1}_i-\bm{p}^{k}_i&=\Delta t \left(  \sum_{j\neq i}\bm{F}_{ij}^{\operatorname{C}}(\bm{q}^{k+1/2})+ (\bm{F}_i^{\operatorname{D}}(\bm{q},\bm{p}))^k \right) +\sigma\sum_{j\neq i}\omega^{\operatorname{R}}(q_{ij}^{k+1/2})\bm{\widehat{q}}_{ij}^{k+1/2}\circ \Delta W_{ij}(t_k),\\
\bm{q}^{k+1}_i-\bm{q}^k_i&=\frac{\Delta t}{m_i}*\frac{\bm{p}_i^k+\bm{p}_i^{k+1}}{2}.
\end{aligned}
\end{equation}

Note that in the latter, $\bm{q}^{k+1/2}$ can also be totally eliminated. We keep it to avoid  heavy arguments for the functions. 

\end{rem}

In the rest of the  paper, we will be focused on the SV-AB schemes \eqref{alg:SV-AB-DPD} (or \eqref{MSV-AB})  for DPD, which include the GW and GCC methods as  special cases. Further studies on the SV-BA schemes and other symplectic methods will be conducted in our future work.
We need only specify the force discretisations $(\bm{F}^{\operatorname{D}}(\bm{q},\bm{p}))^{k_1}$ and $(\bm{F}^{\operatorname{D}}(\bm{q},\bm{p}))^{k_2}$, respectively.  There are certainly many other choices expect for what we introduce below. 

 To recover the conventional GW and GCC methods,  the approximation $\overline{\Delta} W_{ij}(t_k) \approx \overline{\Delta} W_{ij}(t_{k+1/2})$  will have to be employed, and hence
\begin{equation}
\overline{\Delta} W_{ij} \approx \Delta W_{ij}/2,
\end{equation}
 as $\overline{\Delta} W_{ij}(t_k)+\overline{\Delta} W_{ij}(t_{k+1/2})=\Delta W_{ij}(t_k)$. However, it should be noted that this approximation  will change the nature of the schemes in the sense that the increments $\overline{\Delta} W_{ij}(t_k)$ and $\overline{\Delta} W_{ij}(t_{k+1/2})$ are not longer independent; in fact, this approximation is not necessary in practical applications.

\begin{itemize}
\item  {SV-AB-1} is an implicit scheme by choosing 
\begin{equation}
(\bm{F}_i^{\operatorname{D}}(\bm{q},\bm{p}))^{k_1}=\bm{F}^{\operatorname{D}}_i(\bm{q}^k,\bm{p}^k),\quad (\bm{F}^{\operatorname{D}}_i(\bm{q},\bm{p}))^{k_2}=\bm{F}^{\operatorname{D}}_i(\bm{q}^{k+1},\bm{p}^{k+1}).
\end{equation}
For DPD, the dissipative force $\bm{F}^{\operatorname{D}}$ is linear in $\bm{p}$, so the scheme can be written explicitly, in principle. However, one may need to solve a linear system with a sparse coefficient matrix.
\item {SV-AB-2} is an explicit scheme by defining 
\begin{equation}
(\bm{F}_i^{\operatorname{D}}(\bm{q},\bm{p}))^{k_1}=\bm{F}_i^{\operatorname{D}}(\bm{q}^k,\bm{p}^k),\quad (\bm{F}_i^{\operatorname{D}}(\bm{q},\bm{p}))^{k_2}=\bm{F}^{\operatorname{D}}_i(\bm{q}^{k+1},\bm{p}^{k+\lambda}),
\end{equation}
where $\bm{p}^{k+\lambda}$  ($\lambda\in[0,1]$) is defined by 
\begin{equation}\label{GCClambda}
\frac{\bm{p}^{k+\lambda}_i-\bm{p}^k_i}{\Delta t}=\lambda\left[ \sum_{j\neq i}\bm{F}_{ij}^{\operatorname{C}}(\bm{q}^{k})+\bm{F}^{\operatorname{D}}_i(\bm{q}^{k},\bm{p}^k)+\sigma\sum_{j\neq i}\omega^{\operatorname{R}}(q_{ij}^{k})\bm{\widehat{q}}_{ij}^{k}\circ \frac{\Delta W_{ij}(t_k)}{\Delta t}\right].
\end{equation}
This is exactly the GCC method \cite{GCC1999}.
\item  {SV-AB-3} is explicit by specifying 
\begin{equation}
(\bm{F}^{\operatorname{D}}_i(\bm{q},\bm{p}))^{k_1}=\bm{F}^{\operatorname{D}}_i(\bm{q}^k,\bm{p}^{k-1+\lambda}),\quad (\bm{F}^{\operatorname{D}}_i(\bm{q},\bm{p}))^{k_2}=\bm{F}^{\operatorname{D}}_i(\bm{q}^{k+1},\bm{p}^{k+\lambda}),
\end{equation}
where $\bm{p}^{k+\lambda}$  ($\lambda\in[0,1]$) is defined by 
\begin{equation}\label{GWlambda}
\frac{\bm{p}^{k+\lambda}_i-\bm{p}^k_i}{\Delta t}=\lambda\left[ \sum_{j\neq i}\bm{F}_{ij}^{\operatorname{C}}(\bm{q}^{k})+\bm{F}^{\operatorname{D}}_i(\bm{q}^{k},\bm{p}^{k-1+\lambda})+\sigma\sum_{j\neq i}\omega^{\operatorname{R}}(q_{ij}^{k})\bm{\widehat{q}}_{ij}^{k}\circ \frac{\Delta W_{ij}(t_k)}{\Delta t}\right]
\end{equation}
 and the initial value of $\bm{p}^{k+\lambda}$ is $\bm{p}^{\lambda}=\bm{p}^1$ when $k=1$. This is exactly the GW method \cite{Groot1997}.
\item {SV-AB-4} is a generalisation  of the three methods above, which can, in principle,  be expressed explicitly for the DPD: 
\begin{equation}
\begin{aligned}
(\bm{F}^{\operatorname{D}}_i(\bm{q},\bm{p}))^{k_1}&=\bm{F}^{\operatorname{D}}_i(\bm{q}^k,\alpha \bm{p}^{k}+(1-\alpha)\bm{p}^{k-1+\lambda}),\quad \alpha\in[0,1],\\
 (\bm{F}^{\operatorname{D}}_i(\bm{q},\bm{p}))^{k_2}&=\bm{F}^{\operatorname{D}}_i(\bm{q}^{k+1},\beta \bm{p}^{k+\lambda}+(1-\beta)\bm{p}^{k+1}),\quad \beta\in[0,1],
 \end{aligned}
\end{equation}
where $\bm{p}^{k+\lambda}$ ($\lambda\in[0,1]$) is defined by 
\begin{equation}\label{sv4lambda}
\begin{aligned}
\frac{\bm{p}^{k+\lambda}_i-\bm{p}^k_i}{\Delta t} & =\lambda\left[ \sum_{j\neq i}\bm{F}_{ij}^{\operatorname{C}}(\bm{q}^{k})+\bm{F}^{\operatorname{D}}_i(\bm{q}^k,\alpha \bm{p}^{k}+(1-\alpha)\bm{p}^{k-1+\lambda})\right.\\
& \quad \quad \quad \quad \quad  \quad \quad \quad  \left . +~ \sigma\sum_{j\neq i}\omega^{\operatorname{R}}(q_{ij}^{k})\bm{\widehat{q}}_{ij}^{k}\circ \frac{\Delta W_{ij}(t_k)}{{\Delta t}}\right].
\end{aligned}
\end{equation}
It reduces to the SV-AB-1 method for $\alpha=1,\beta=0$, to the SV-AB-2 (GCC) method for $\alpha=1,\beta=1$ and to the SV-AB-3 (GW) method for $\alpha=0,\beta=1$.
\item  {SV-AB-5} is explicit by choosing 
\begin{equation}
\begin{aligned}
(\bm{F}_i^{\operatorname{D}}(\bm{q},\bm{p}))^{k_1}&=\bm{F}_i^{\operatorname{D}}(\bm{q}^k,\bm{p}^{k-1+\lambda_1}),\quad \lambda_1\in[0,1],\\
 (\bm{F}_i^{\operatorname{D}}(\bm{q},\bm{p}))^{k_2}&=\bm{F}_i^{\operatorname{D}}(\bm{q}^{k+1},\bm{p}^{k+\lambda_2}),\quad \lambda_2\in[0,1],
 \end{aligned}
\end{equation}
where 
\begin{equation}\label{lambda12}
\begin{aligned}
\frac{\bm{p}^{k+\lambda_1}_i-\bm{p}^k_i}{\Delta t}&=\lambda_1\left[ \sum_{j\neq i}\bm{F}_{ij}^{\operatorname{C}}(\bm{q}^{k})+\bm{F}^{\operatorname{D}}_i(\bm{q}^k,\bm{p}^{k-1+\lambda_1})+\sigma\sum_{j\neq i}\omega^{\operatorname{R}}(q_{ij}^{k})\bm{\widehat{q}}_{ij}^{k}\circ \frac{\Delta W_{ij}(t_k)}{\Delta t}\right],\\
\frac{\bm{p}^{k+\lambda_2}_i-\bm{p}^{k}_i}{\Delta t}&=\lambda_2\left[ \sum_{j\neq i}\bm{F}_{ij}^{\operatorname{C}}(\bm{q}^{k})+\bm{F}^{\operatorname{D}}_i(\bm{q}^{k},\bm{p}^{k-1+\lambda_1})+\sigma\sum_{j\neq i}\omega^{\operatorname{R}}(q_{ij}^{k})\bm{\widehat{q}}_{ij}^{k}\circ \frac{\Delta W_{ij}(t_k)}{\Delta t}\right].
\end{aligned}
\end{equation} 
When $\lambda_1=\lambda_2$, it reduces to the SV-AB-3 (GW) method. 
\item {SV-AB-6} is a simultaneous generalisation of {SV-AB-4} and {SV-AB-5}, which can be written in an explicit form for the DPD:  
\begin{equation}
\begin{aligned}
(\bm{F}^{\operatorname{D}}_i(\bm{q},\bm{p}))^{k_1}&=\bm{F}^{\operatorname{D}}_i(\bm{q}^k,\alpha \bm{p}^{k}+(1-\alpha)\bm{p}^{k-1+\lambda_1}),\quad \alpha\in[0,1],\lambda_1\in[0,1],\\
 (\bm{F}^{\operatorname{D}}_i(\bm{q},\bm{p}))^{k_2}&=\bm{F}^{\operatorname{D}}_i(\bm{q}^{k+1},\beta \bm{p}^{k+\lambda_2}+(1-\beta)\bm{p}^{k+1}),\quad \beta\in[0,1],\lambda_2\in[0,1],
 \end{aligned}
\end{equation}
where $\bm{q}^{k+\lambda_1}$ and $\bm{q}^{k+\lambda_2}$ are given by
 \begin{equation}
\begin{aligned}
\frac{\bm{p}^{k+\lambda_1}_i-\bm{p}^k_i}{\Delta t}&=\lambda_1\left[ \sum_{j\neq i}\bm{F}_{ij}^{\operatorname{C}}(\bm{q}^{k})+\bm{F}^{\operatorname{D}}_i(\bm{q}^k,\alpha \bm{p}^{k}+(1-\alpha)\bm{p}^{k-1+\lambda_1})\right.\\
&\quad\quad\quad\quad\quad\quad\quad \quad \left.+ ~\sigma\sum_{j\neq i}\omega^{\operatorname{R}}(q_{ij}^{k})\bm{\widehat{q}}_{ij}^{k}\circ \frac{\Delta W_{ij}(t_k)}{\Delta t}\right],\\
\frac{\bm{p}^{k+\lambda_2}_i-\bm{p}^{k}_i}{\Delta t}&=\lambda_2\left[ \sum_{j\neq i}\bm{F}_{ij}^{\operatorname{C}}(\bm{q}^{k})+\bm{F}^{\operatorname{D}}_i(\bm{q}^{k},\alpha \bm{p}^{k}+(1-\alpha)\bm{p}^{k-1+\lambda_1})\right.\\
&\quad\quad\quad\quad\quad\quad\quad \quad\left.+~\sigma\sum_{j\neq i}\omega^{\operatorname{R}}(q_{ij}^{k})\bm{\widehat{q}}_{ij}^{k}\circ \frac{\Delta W_{ij}(t_k)}{\Delta t}\right].
\end{aligned}
\end{equation} 
When $\lambda_1=\lambda_2$, it becomes SV-AB-4, while when $\alpha=0,\beta=1$, it becomes SV-AB-5.
\end{itemize}

\begin{rem}\label{rem:SVss}
 The schemes {SV-AB-1}$\sim${6} are related through the following diagram:
\begin{figure}[H]
\centering
\begin{minipage}{10cm}
\setlength{\unitlength}{1cm}
\begin{picture}(18,7.5)
\put(-0.8,7){{SV-AB-1}}
\put(3.4,7){{SV-AB-2} (GCC)}
\put(8.0,7){{SV-AB-3} (GW)}
\put(4.2,4.5){{SV-AB-4}}
\put(8.6,4.5){{SV-AB-5}}
\put(6.4,2){{SV-AB-6}}

\put(5,5){\vector(0,1){1.8}}
\put(5,5){\vector(4,1.7){4.4}}
\put(5,5){\vector(-5.5,2){4.9}}

\put(7.2,2.5){\vector(-2,1.7){2.1}}
\put(7.2,2.5){\vector(2,1.6){2.3}}

\put(9.45,5){\vector(0,1){1.8}}

\put(5.1,6.2){$\alpha=1$}
\put(5.1,5.7){$\beta=1$}

\put(1.2,6.2){$\alpha=1$}
\put(1.2,5.7){$\beta=0$}

\put(7,6.2){$\alpha=0$}
\put(7,5.7){$\beta=1$}

\put(9.53,5.9){$\lambda_1=\lambda_2$}

\put(4.7,3.2){$\lambda_1=\lambda_2$}

\put(7.7,3.6){$\alpha=0$}
\put(7.7,3.1){$\beta=1$}

\end{picture}
\end{minipage}
\vspace{-1.5cm}
\caption{Relations of the schemes {SV-AB-1}$\sim${6}.} 
\label{fig:svab}
\end{figure}
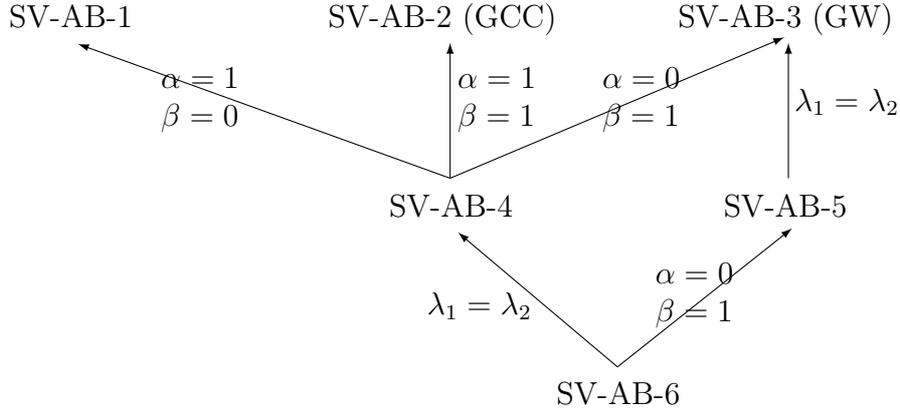

\end{rem}

We summarize some general features of the schemes in Fig. \ref{fig:svab} as follows.
\begin{itemize}
\item Explicit schemes: {SV-AB-2} (GCC), {SV-AB-3} (GW), {SV-AB-5}. The other schemes are implicit but can be written explicitly for the DPD by solving a sparse linear system.
\item Number of independent parameters: {SV-AB-1}: 0, {SV-AB-2} (GCC): 1, {SV-AB-3} (GW): 1, {SV-AB-4}: 3, {SV-AB-5}: 2,  {SV-AB-6}: 4. 
\end{itemize} 

\subsection{Long-time behaviour of the SV-AB methods}
When no external forces are involved, it was noticed that the Euler-A method \eqref{alg:EA} and the Euler-B method \eqref{alg:EB} are not convergent in the mean-square sense when the Hamiltonian functions $h_{ij}=h_{ij}(\bm{q},\bm{p})$ depend on both the positions and momenta \cite{HT2018}. 
If we further assume that $h_{ij}=h_{ij}(\bm{q})$ only depend on the positions, the Euler-A and Euler-B methods are both convergent and hence are the SV methods. 

In this subsection, we will numerically show the convergence of the {SV-AB}-1$\sim$6 methods by studying the damped Kubo oscillator, a stochastic Hamiltonian system whose Hamiltonians are separable given by 
\begin{equation}
H(q,p)=\frac{p^2}{2}+\frac{q^2}{2},\quad h(q,p)=\sigma\left(\frac{p^2}{2}+\frac{q^2}{2}\right).
\end{equation}
Here $\sigma$ is the noise intensity. As its solution can be calculated analytically, it has often been used for the validation of numerical methods (e.g., \cite{HT2018,MRT2002}). By employing the forces 
\begin{equation}
F^{\operatorname{D}}=-\varepsilon p, \quad F^{\operatorname{SD}}=-\varepsilon\sigma p
\end{equation}
with $\varepsilon$ the nonnegative damping coefficient, the damped Kubo oscillator has the following exact solution \cite{KT2021}
\begin{equation}
\begin{aligned}
\overline{q}(t)&=q_0\exp\left(-\frac{\varepsilon}{2}(t+\sigma W(t))\right)\cos\omega\left(t+\sigma W(t)\right)\\
&\quad  +\frac{1}{\omega}\left(p_0+\frac{\varepsilon}{2}q_0\right)\exp\left(-\frac{\varepsilon}{2}(t+\sigma W(t))\right)\sin\omega\left(t+\sigma W(t)\right),\\
\overline{p}(t)&=p_0\exp\left(-\frac{\varepsilon}{2}(t+\sigma W(t))\right)\cos\omega\left(t+\sigma W(t)\right)\\
&\quad  -\frac{1}{\omega}\left(q_0+\frac{\varepsilon}{2}p_0\right)\exp\left(-\frac{\varepsilon}{2}(t+\sigma W(t))\right)\sin\omega\left(t+\sigma W(t)\right),
\end{aligned}
\end{equation}
where $(q_0,p_0)$ are the initial conditions, the angular frequency is $\omega={\sqrt{4-\varepsilon^2}}/{2}$ by assuming $\varepsilon<2$. The expected value of the Hamiltonian $H$ is given by
\begin{equation}\label{exhal}
\begin{aligned}
E(H(\overline{q}(t),\overline{p}(t)))&=a\exp\left(-\frac{\varepsilon(2-\varepsilon\sigma^2)t}{2}\right)\\
&\quad +\exp\left(-\left((2-\varepsilon^2)\sigma^2+\varepsilon\right)t\right)\left(b\cos\left(2(1-\varepsilon\sigma^2)\omega t\right)+c\sin\left(2(1-\varepsilon\sigma^2)\omega t\right)\right),
\end{aligned}
\end{equation}
where 
\begin{equation}
a=\frac{2(q_0^2+p_0^2+\varepsilon q_0p_0)}{4-\varepsilon^2},\quad b=-\frac{\varepsilon^2(q_0^2+p_0^2)+4\varepsilon q_0p_0}{2(4-\varepsilon^2)},\quad c=\frac{\varepsilon(q_0^2-p_0^2)}{2\sqrt{4-\varepsilon^2}}.
\end{equation}

\begin{figure}[tbhp]
\centering
\includegraphics[width=.9\linewidth]{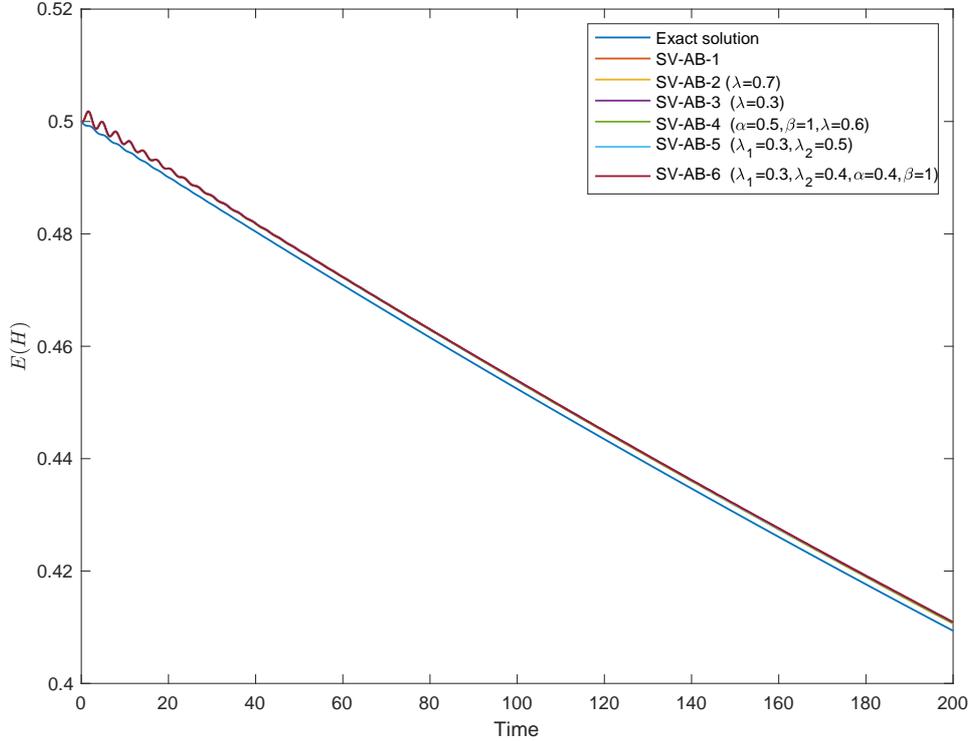}
\caption{Mean Hamiltonians of the SV-AB-1, SV-AB-2 ($\lambda=0.7$), SV-AB-3 ($\lambda=0.3$), SV-AB-4 ($\alpha=0.5,\beta=1,\lambda=0.6$), SV-AB-5 ($\lambda_1=0.3,\lambda_2=0.5$), and SV-AB-6 ($\lambda_1=0.3,\lambda_2=0.4,\alpha=0.4,\beta=1$) methods. To clearly show the tendency of the time evolution and the fluctuating behaviour, the first 200 seconds are plotted here.}
\label{fig:kuboh}
\end{figure}

\begin{figure}[tbhp]
\centering
\includegraphics[width=.9\linewidth]{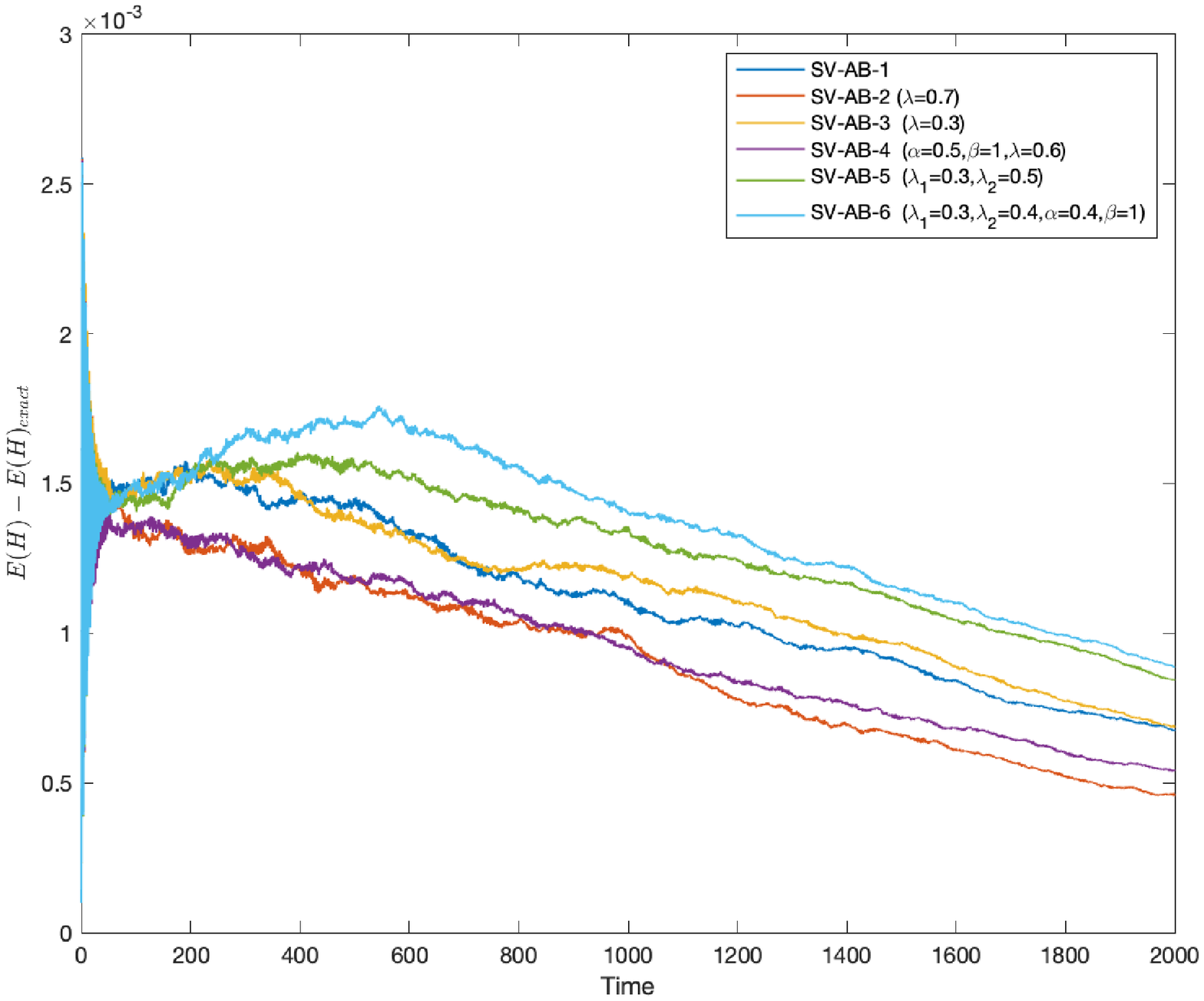}
\caption{The difference between numerical mean Hamiltonians and the exact Hamiltonian for the SV-AB-1, SV-AB-2 ($\lambda=0.7$), SV-AB-3 ($\lambda=0.3$), SV-AB-4 ($\alpha=0.5,\beta=1,\lambda=0.6$), SV-AB-5 ($\lambda_1=0.3,\lambda_2=0.5$), and SV-AB-6 ($\lambda_1=0.3,\lambda_2=0.4,\alpha=0.4,\beta=1$) methods.}
\label{fig:kubohe}
\end{figure}

In the simulations, the initial conditions are $q_0=0,p_0=1$, the noise intensity is $\sigma=0.2$ and the damping coefficient is $\varepsilon=0.001$.  Time step is $\Delta t=0.1$ for a time span $[0,2000]$. 
For simplicity, discretisation of $F^{\operatorname{SD}}(q,p)$ is chosen as $F^{\operatorname{SD}}(q^k,p^k)$ at each step $k$ for all numerical methods. Furthermore, we pick one special choice of the parameters $\alpha,\beta,\lambda$ for each method as shown in the figures and in each case $2,000$ sample paths are generated. Figs. \ref{fig:kuboh} and \ref{fig:kubohe} show the mean Hamiltonians of the SV-AB methods and their differences with respect to the exact Hamiltonian \eqref{exhal}.  Fluctuating behaviour of the energy can be noticed. In particular, Fig. \ref{fig:kubohe} shows that order of the error is approximately $10^{-3}$ and it tends to become smaller in a long time after a relatively  stronger vibration at the beginning.  


\section{Applications to DPD simulations}
\label{sec:app}

Although all the SV-AB schemes proposed above can be made explicit for the DPD, further efforts may be needed to achieve the corresponding explicit representatives, in particular, by solving a huge sparse linear system.  For simplicity, we will be focused on the explicit  {SV-AB-2} (GCC) and   {SV-AB-4} ($\beta=1$) methods in comparison with the {SV-AB-3} (GW) method. Recall that {SV-AB-4} ($\beta=1$) reduces to the {SV-AB-3} (GW) with $\alpha=0$ and reduces to SV-AB-2 (GCC) with  $\alpha=1$ (see Fig. \ref{fig:svab}).
In our simulations, the total number of fluid particles of the same mass $m$ is set to $3,000$ with  $a=25 k_\mathrm{B}T^*$, where $a$ is the repulsive parameter (i.e., $a_{ij}=a$ for all $i\neq j$) to determine the magnitude of the conservative force $\bm{F}^{\operatorname{C}}$, $T^*$ is the set temperature and $k_\mathrm{B}$ is the Boltzmann constant.
The noise parameter $\sigma$ and the friction parameter $\gamma$ are set to $3.0$ and $4.5$, respectively. 
All simulations are performed under the condition of constant-volume and constant-temperature, i.e., the canonical ensemble is generated.
The size of simulation box is 10 $\times 10 \times 10 q_{\mathrm{c}}^3$. 
The periodic boundary condition is applied in all three dimensions. 
 Here, $q_{\mathrm{c}}$ is the cutoff distance, which is the unit length in the DPD simulation.
The initial configuration is random, and the initial momentum is set appropriately so that the temperature would satisfy the Boltzmann distribution for the set temperature satisfying $ k_\mathrm{B}T^*=1.0$. This gives the repulsion parameter $a=25$, yielding the compressibility of water.
Although Groot and Warren reported that there was no statistical difference between simulations using uniform random numbers and those using Gaussian random numbers \cite{Groot1997}, we use a Gaussian distribution to generate the random numbers in the current simulations.

We examined twenty cases with the time step size $\Delta t$ ranging from $0.001$ to $0.16 \tau$. Here, we use reduced units for the cutoff radius $q_{\mathrm{c}}$, the particle mass $m$, and the energy $k_\mathrm{B}T$.  
Hence, the time unit is defined as $\tau = \sqrt{m q_{\mathrm{c}}^2/k_\mathrm{B}T}$.
All cases were simulated for at least $1,000 \tau$, and the last $16 \%$ were used as statistical data. 
Note that we were not able to calculate exactly $1,000 \tau$ for all $\Delta t$ and the first 84 \% of the data was discarded to equilibrate the system sufficiently.
As a comparison of the accuracy of the formulations, the kinetic temperature $k_\mathrm{B}T = \left\langle \bm{v}^2 \right\rangle/3$ was calculated  and its difference from the set temperature  $ k_\mathrm{B}T^*=1.0$ was examined, where $\langle \cdot \rangle $ is  the average over all particles in the simulations and $\bm{v}=\bm{p}/m$.
 Since the simulation was performed with a canonical ensemble, temperature of the system will fluctuate around a certain average value after reaching the equilibrium state. In the simulations, the average value is the set temperature, which satisfies  $k_\mathrm{B}T^*=1.0$. 

\begin{figure}[tbhp]
\centering
\includegraphics[width=.8\linewidth]{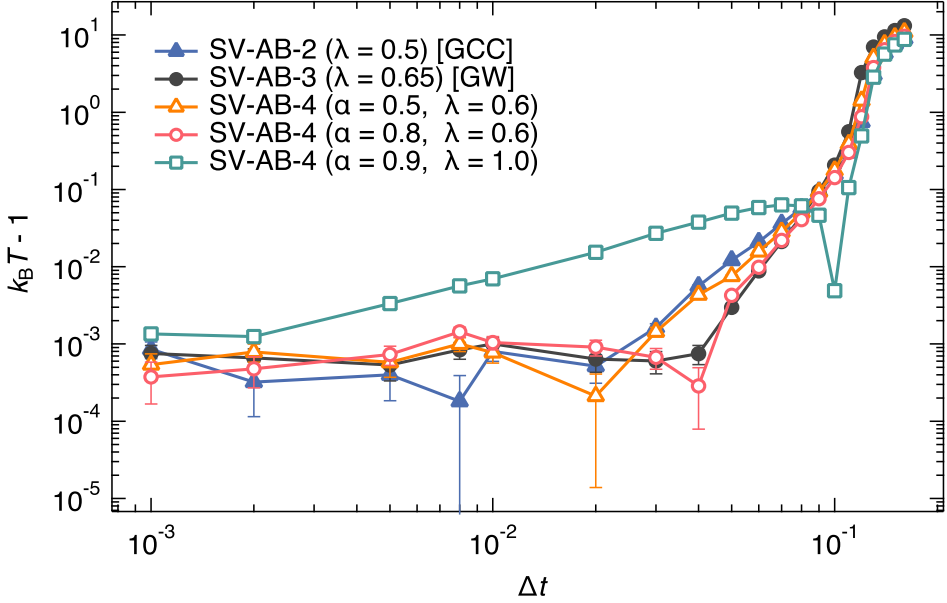}
\caption{Kinetic temperature versus time step. Curves represent representative results for the {SV-AB-2 (GCC)}, {SV-AB-3 (GW)}, and {SV-AB-4} ($\beta=1$) schemes.  Note that the kinetic temperature is averaged over time after equilibration.}
\label{fig02}
\end{figure}

Fig. \ref{fig02} plots the artificial kinetic temperature increase of the {SV-AB-2 (GCC)}, {SV-AB-3 (GW)}, and {SV-AB-4} ($\beta=1$) schemes with representative parameters. 
For results for all parameters, please refer to Figs. S1--S3 in Supporting Information.
It is confirmed that the statistical error of the temperature is less than $1 \%$, i.e., $k_\mathrm{B}T - 1 < 10^{-2}$, for all schemes when $\Delta t$ is less than $0.01$. 

Let us firstly compare the {SV-AB-3 (GW)} scheme with the {SV-AB-2 (GCC)} scheme. We consider that $\lambda = 0.5$ and $0.65$ are the representative parameters of the {SV-AB-2 (GCC)} and {SV-AB-3 (GW)} schemes, respectively.  When $\Delta t< 2\times 10^{-2}$,  in several cases the error of {SV-AB-2 (GCC)} ($\lambda=0.5$) is smaller than that of the  {SV-AB-3 (GW)}  with $\lambda=0.65$. When $\Delta t >3\times10^{-2}$, error of {SV-AB-2 (GCC)} ($\lambda=0.5$) jumps to bigger than $0.1\%$. However, when the time step becomes ever bigger, for instance $\Delta t>10^{-1}$, error of {SV-AB-2 (GCC)} ($\lambda=0.5$) is smaller; one should be noted that  error for these cases is  probably too big for practical simulations.

Now consider the {SV-AB-4} ($\beta=1$) scheme. 
For all $\alpha$s, the error tends to be the smallest around $\lambda=0.6$.
When $\Delta t < 10^{-2}$, the error of the {SV-AB-4} ($\beta=1$) is similar to that of the {SV-AB-3 (GW)}  ($\lambda=0.65$). 
As $\Delta t$ increases, the error also increases. The accuracy of schemes with $\alpha=0.8$ and $\alpha=0.5$ interchanges at some point  as $\Delta t$  increases: for smaller $\Delta t$, the error of $\alpha = 0.5$ case is smaller, while for larger $\Delta t$, the error of $\alpha = 0.8$ case becomes smaller. The maximum $\Delta t$ that shows an accuracy of less than $1 \%$ error is $0.06$, which is the same as the {SV-AB-3 (GW)} ($\lambda=0.65$), but for  $\Delta t = 0.04$, its accuracy is higher than the {SV-AB-3 (GW)}  ($\lambda=0.65$).
On the other hand, when $\Delta t < 7 \times 10^{-2}$ and $\lambda = 1.0$, the error of {SV-AB-4} ($\beta=1$) is larger than that of the {SV-AB-3 (GW)}  ($\lambda=0.65$) for all $\alpha$s. However, when $\Delta t = 0.1$, the error is approximately $0.5 \%$, which is highly accurate. 
Unfortunately, further studies are needed before this can be applied in practical simulations easily. 

\begin{figure}[tbhp]
\centering
\includegraphics[width=.8\linewidth]{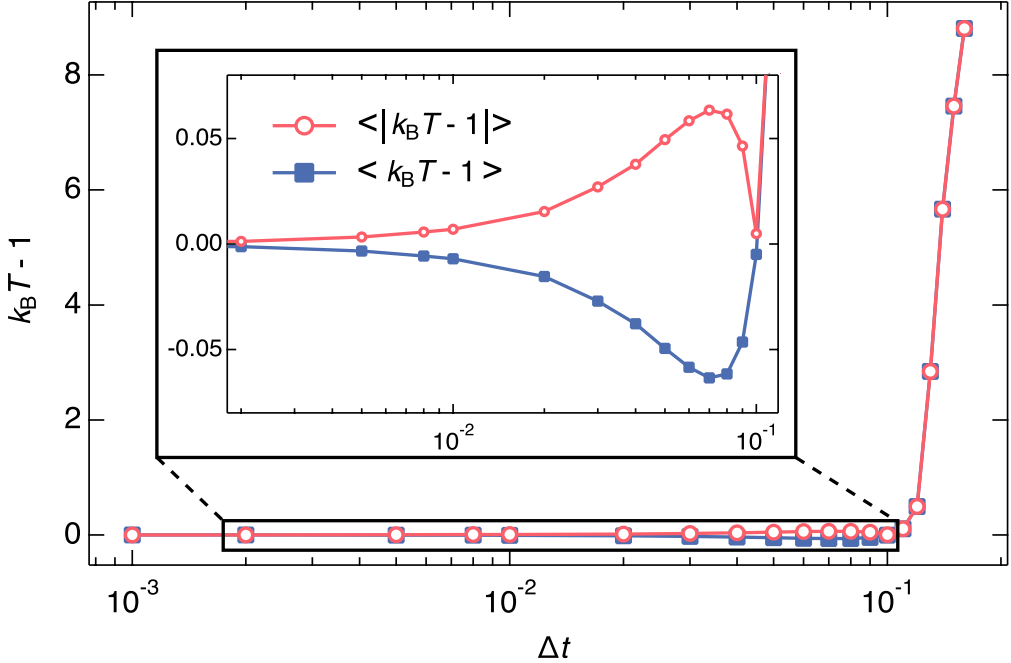}
\caption{Kinetic temperature versus time step on linear $y$-axis. {SV-AB-4} ($\beta=1$) scheme with $\alpha = 0.9$ and $\lambda = 1.0$. Red and blue curves represent statical error of kinetic temperatures before and after taking absolute values.  Note that the kinetic temperature is averaged over time after equilibration.}
\label{fig03}
\end{figure}

Simulations of the {SV-AB-4} ($\beta=1$) for $\alpha = 0.9$ and $\lambda = 1.0$ 
are shown in Fig. \ref{fig03} with the vertical axis illustrated in linear scale. 
Blue and red curves show the error and the absolute error respectively. 
Note that the green curve in Fig. \ref{fig02} and the red curve in Fig. \ref{fig03} coincide.
As shown in the zoomed-up view inserted in Fig. \ref{fig03}, starting  at around $\Delta t=2\times10^{-2}$, the statical error becomes bigger and bigger in negative values as $\Delta t$ increases.  It attains $-0.05$ at $\Delta t = 7\times10^{-2}$ and then becomes smaller towards the positive direction as $\Delta t$ increases further. Finally at $\Delta t = 0.1$,  the error shifts from negative  to positive and  it is expected that the error is not small for all schemes. This phenomenon that the error is shifting between negative and positive is  observed in all three methods including the GW and the GCC. As a consequence, in practical applications, it is necessary to adopt the value of time step size $\Delta t$ within the permissible error during when error becomes larger as $\Delta t$ becomes larger, instead of the value of $\Delta t$ with the smallest error.

\section{Conclusions and outlook}
\label{sec:con}

In conclusion, we proposed a novel stochastic Hamiltonian formulation (SHF) with matrix noise and subject to external forces, which was found applicable to DPD simulations as the DPD system could be obtained from the corresponding stochastic Lagrange--d'Alembert principle by introducing proper Hamiltonian functions and dissipative forces. In particular, we extended the well-known symplectic SV scheme for conservative Hamiltonian systems to the SHF as composites of the  Euler-A and Euler-B methods. By discretising the dissipative forces properly, several simple families of SV methods were constructed and especially the SV-AB methods were focused which were derived as the composite $\text{(Euler-A)}\circ \text{(Euler-B)}$. By studying the damped Kubo oscillator, the fluctuating behaviour and damping energy/Hamiltonian dissipation were realised with order of error approximately $10^{-3}$ between the numerical  and exact Hamiltonians.
For DPD simulations, the SV-AB methods include the conventional GW and GCC methods as special cases. Simulations of a novel two-parametric explicit schemes were conducted and compared with the GW and  GCC methods.  As time step varies, some of the novel schemes were advantageous over the GW method but unfortunately no global advantage was realised. It was also observed  that for all schemes as the time step increases the error can shift between positive and negative values, that requires one to choose a time step in practical  applications more carefully. 
 
Beside the SV methods proposed in the current study, thanks to the SHF a variety of other effective  structure-preserving methods may be extended as well, for instance, symplectic partitioned Runge--Kutta methods and variational integrators. These are part of our current and future studies including their applications to the DPD and other relevant stochastic  physical systems. From the theoretical viewpoint, it is worthwhile to study further the geometric and algebraic structures of the SHF, for instance, conformal symplectic structures, generating functions, symmetries and Noether's conserved quantities.

\section*{Declaration of competing interest}
The authors declare that they have no known competing financial interests or personal relationships that could have appeared to influence the work reported in this paper.


\section*{Acknowledgements}
This work was partially supported by JSPS KAKENHI Grant Number JP20K14365, JST-CREST Grant Number JPMJCR1914, and Keio Gijuku Fukuzawa Memorial Fund. We thank the anonymous referees for their constructive comments.

\bibliographystyle{elsarticle-num}

\bibliography{article_new}

\end{document}